\documentclass[a4paper,journal]{IEEEtran}
\usepackage{amsmath,amsfonts}
\usepackage{algorithmic}
\usepackage{algorithm}
\usepackage{array}
\usepackage[caption=false,font=normalsize,labelfont=sf,textfont=sf]{subfig}
\usepackage{textcomp}
\usepackage{stfloats}
\usepackage{url}
\usepackage{verbatim}
\usepackage{graphicx}
\usepackage{cite}
\usepackage[numbers, sort&compress, square]{natbib}
\usepackage{enumerate}
\usepackage{accents}
\usepackage{multirow}
\usepackage{amsbsy}
\usepackage{latexsym}
\usepackage{svg}
\usepackage{float}
\usepackage{adjustbox}
\newcommand\scalemath[2]{\scalebox{#1}{\mbox{\ensuremath{\displaystyle #2}}}}
\usepackage{pifont}

\newcommand{\setK}{\mathcal{K}}

\newcommand{\Pgmax}{\overline{P_g}}
\newcommand{\Pgmin}{\underline{P_g}}

\makeatletter
\newcommand{\douwidehat}[2]{%
  \sbox0{$\m@th#1\widehat{\hphantom{#2}}$}%
  \sbox2{$\m@th#1x$}
  \sbox4{$\m@th#1#2$}
  \dimen0=\ht0
  \advance\dimen0 -.8\ht2
  \dimen2=\dp4
  \rlap{%
    \raisebox{\dimexpr\dimen0-\dimen2}{%
      \scalebox{1}[-1]{\box0}%
    }%
  }%
  {#2}%
}
\makeatother

\usepackage[pdfstartview=XYZ,bookmarks=true,colorlinks=true,linkcolor=blue,urlcolor=blue,citecolor=blue,pdftex,bookmarks=false,linktocpage=true, hyperindex=true]{hyperref}
\usepackage{orcidlink}
\begin{document}
\title{AC-DC Power Systems Optimization with Droop Control Smooth Approximation}


\author{G. Mohy-ud-din \orcidlink{0000-0003-2738-9338},{~\IEEEmembership{Member, IEEE,}}
    R. Heidari \orcidlink{0000-0003-4835-2614},{}
    F. Geth \orcidlink{0000-0002-9345-2959},{~\IEEEmembership{Member, IEEE,}} \\
    H. Ergun \orcidlink{0000-0001-5171-1986},{~\IEEEmembership{Senior Member, IEEE}
    and
    S. M. Muslem Uddin \orcidlink{0000-0002-4815-4872} {}
    }

\thanks{G. Mohy-ud-din, R. Heidari, and, S. M. Muslem Uddin are with the Energy Business Unit, Commonwealth Scientific and Industrial Research Organisation (CSIRO), Newcastle, NSW 2300, Australia, (e-mail: ghulam.mohyuddin@csiro.au; rahmat.heidarihaei@csiro.au; smmuslem.uddin@csiro.au).

F. Geth is with  GridQube, Brisbane, QLD 4300, Australia, (e-mail: frederik.geth@gridqube.com).

H. Ergun is with  ELECTA, Department of Electrical Engineering, KU Leuven, Leuven 3000, Belgium, and EnergyVille, Genk, Belgium (e-mail: hakan.ergun@kuleuven.be).
}}



\maketitle

\thispagestyle{empty}

\begin{abstract}
This paper addresses the challenges of embedding common droop control characteristics in ac-dc power system steady-state simulation and optimization problems.
We propose a smooth approximation methodology to construct differentiable functions that encode the attributes of piecewise linear droop control with saturation. We transform the nonsmooth droop curves into smooth nonlinear equality constraints, solvable with Newton methods and interior point solvers. These constraints are then added to power flow, optimal power flow, and security-constrained optimal power flow problems in ac-dc power systems. The results demonstrate significant improvements in accuracy in terms of power sharing response, voltage regulation, and system efficiency, while outperforming existing mixed-integer formulations in computational efficiency.
\end{abstract}



\section{Introduction}

Integrating renewable energy into power systems presents significant challenges for voltage control, which is essential for maintaining network security \cite{voltage_control}. In ac systems, voltage adjustments are made through reactive power control, while dc systems rely on active power control via HVDC converters. Due to the inherent limitations in transmitting reactive power in ac systems, typical voltage regulation involves various compensating devices: passive compensators like shunt reactors and capacitors, active compensators such as synchronous condensers, static var compensators, FACTS devices, and tap-changing transformers. Advanced control strategies such as droop control are crucial as they enable power-sharing and voltage regulation without relying on communication, thus facilitating distributed and reliable system control \cite{droopControl}.

Droop control is not limited to active power and frequency regulation (P--f droop); it also includes reactive power-voltage control (Q--V droop), commonly employed in synchronous generators and commercial-scale inverters \cite{fpqv_gdroop}, \cite{fpqv_idroop}. 
By adjusting the generator or inverter output based on deviations in voltage and frequency, these control mechanisms help stabilize the grid by balancing reactive power supply and demand. Typically, grid-following inverters operate with PQ control, and grid-forming inverters are droop-controlled \cite{PQandDroop}. However, different technologies and configurations employ more advanced strategies such as adaptive, coordinated consensus, current-limiting, and universal droop control \cite{adaptive}, \cite{consensus}, \cite{currentlimiting}, \cite{universal}. Similarly, active power-voltage (P--Vdc) and reactive power voltage (Q--Vac) droop controls are implemented in HVDC converters to regulate voltage levels by adjusting the active power output of HVDC converters in response to dc-side voltage variations, and reactive power injection in response to ac-side voltage variations thus ensuring stability during transient conditions \cite{vsc_droop}, \cite{lcc_droop}.
The two primary converter technologies, Voltage Source Converter (VSC) and Line Commutated Converter (LCC), have different characteristics. Therefore, in mixed configurations such as parallel VSC-LCC \cite{parallel_lcc_vsc}, cascaded VSC-LCC \cite{cascaded_lcc_vsc}, and multi-infeed VSC-LCC-based dc power systems \cite{multi-infeed_lcc_vsc}, coordinated and unified droop control strategies are implemented. 


The key to the performance of droop control is the  characteristic function. Traditionally, linear droop characteristics are implemented, however, due to a constant slope, it faces a trade-off between strict voltage regulation and accurate power-sharing. A steeper slope improves power-sharing but increases voltage deviation under high loads, while a gentler slope enhances voltage regulation but reduces power-sharing effectiveness. An optimization to determine the droop slope ensuring operation security against credible contingency events is an effective approach addressing this trade-off~\cite{ergunRobust}. 

Various nonlinear droop functions such as second-order (parabola, inverse parabola, ellipse) and higher-order polynomials have been explored. With these polynomial functions, the slope of the curves can be precisely calibrated from no-load to high-load conditions by tuning the coefficients \cite{droopControl}, \cite{investigationNonlinearDroop}. They enable more precise fitting of droop gains across multiple curves improving voltage regulation but sacrificing power-sharing at different parts of the curves and vice-versa.

\emph{Piecewise} linear droop functions address the previous challenges by offering customizable tuning and ideal droop characteristics by combining the strengths of different droop curves \cite{investigationNonlinearDroop}. It divides the droop curve into multiple segments with adjustable slopes and supports deadbands. These deadbands are essential for maintaining stability, reducing control chatter, and improving system response by introducing tolerance around the set-point. While piecewise linear droop control is effective for various power systems, abrupt slope changes between segments can be problematic during load transitions, thus, a smooth droop curve is preferred \cite{investigationNonlinearDroop}.

\subsection{Challenges of Modelling Droop Controls in Optimal Power Flow Problems}

Piecewise linear droop controls are difficult to implement for several reasons. The non-differentiable  nature of these functions poses challenges for control system design, tuning, and optimization. The control system design requires a solution to the complete representation of the power system including steady-state operation and dynamic controls formulated as a system of differential and algebraic equations, where non-differentiable functions add further complexity \cite{DAEdroopApprox1}, \cite{DAEdroopApprox2}. 

The implementation of piecewise linear droop control in steady-state nonlinear power systems problems, such as power flow (PF) \cite{PF}, optimal power flow (OPF) \cite{pmacdc}, and security-constrained optimal power flow (SCOPF) \cite {SCOPF}  is  challenging. The PF problems typically use numerical solvers which only support twice differentiable functions in the model. In OPF and SCOPF problems, the typical method to encode piecewise linear droop control depends on integer variables. 
The inclusion of integer variables makes the overall problem with the ac/dc physics a mixed-integer nonlinear programming (MINLP) problem and therefore it is much more challenging to find optimal solutions \cite{computationalIssues}. 

\vspace{-0.5cm}

\subsection{Scope and Contributions}
To address these challenges, we propose a methodology to convert the piece-wise linear droop characteristics to smooth (differentiable everywhere)  functions.  We use a technique inspired by the machine learning literature, which entails encoding the piecewise functions through the \emph{ReLU}\footnote{`rectified linear unit'} function. After encoding the constraints in that fashion, we  apply a substitution with a function that smoothly approximates the ReLU, i.e. the \emph{softplus} function with a tuneable sharpness setting. 

Overall, this process converts the nonsmooth  droop characteristics into smooth but approximate constraints with tuneable approximation error, which can be efficiently solved using  Newton's method or derivative-based interior point solvers. Formulations are provided for different droop controls in ac-dc power systems and implemented in PF, OPF, and SCOPF problems to demonstrate the computational efficiency and accuracy in comparison to the mixed-integer formulations. 
%
%
In summary, the key contributions are:
\begin{itemize}
    \item A technique to derive smooth models encoding the attributes of piece-wise linear droop control including saturation. This enables more realistic SCOPF models that ensure better power-sharing, voltage regulation, and, system efficiency without expensive computation.
    \item The implementation is made available, with support for different droop controls in ac-dc power systems in PF, OPF, and SCOPF problems demonstrating computation efficiency and accuracy in comparison to the mixed-integer formulations.
\end{itemize}

This paper is organized as follows. Section II describes the proposed smooth droop control functions. Sections III and IV present the simulation results and conclusions.

\begin{figure}[!t]
    \centering
    \includegraphics[width=2.9in]{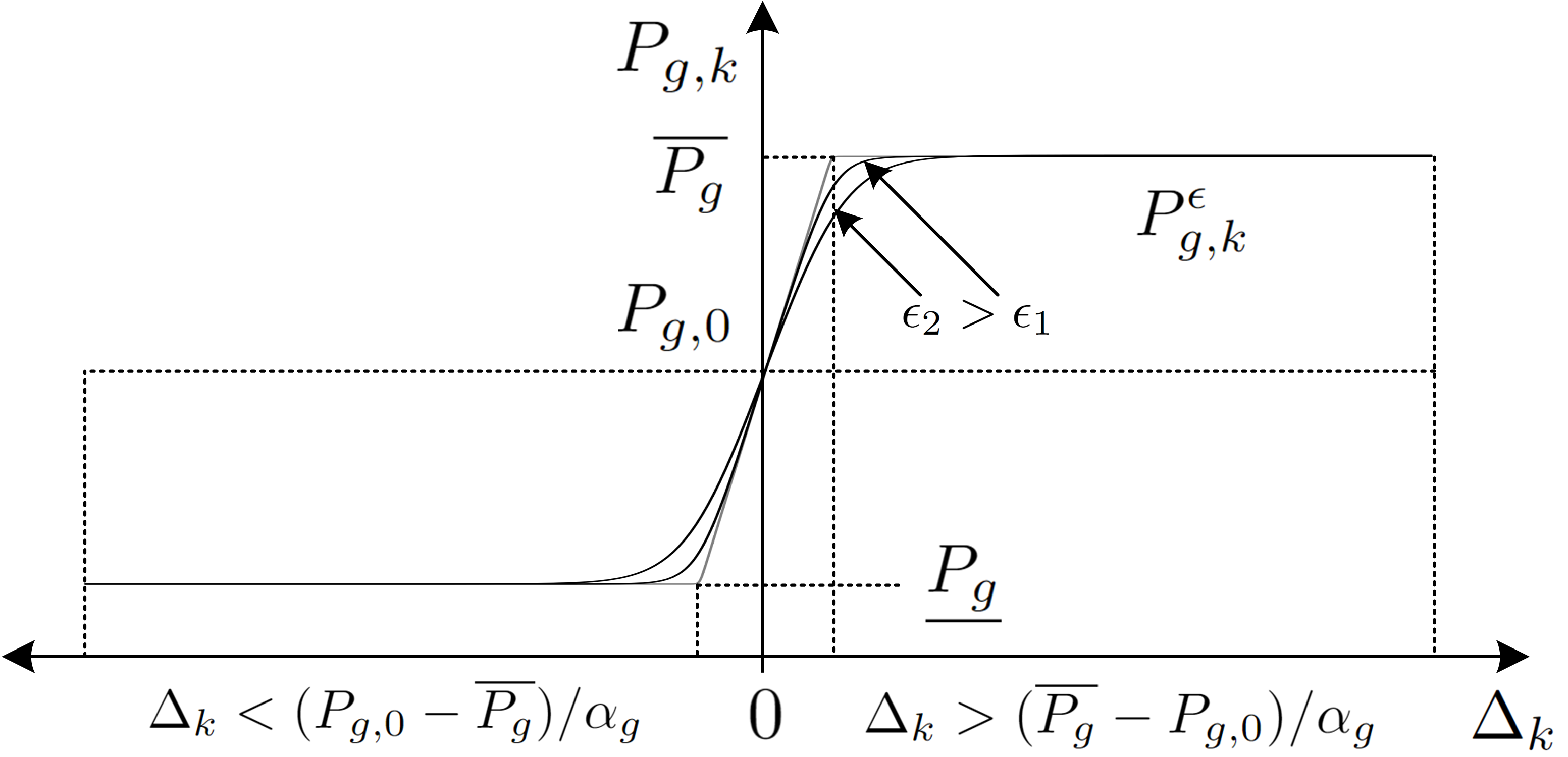}
    \caption{Original and approximate generator active power response functions.}
    \label{fig:genP}
\end{figure}
\begin{figure}[!t]
    \centering
    \includegraphics[width=2.9in]{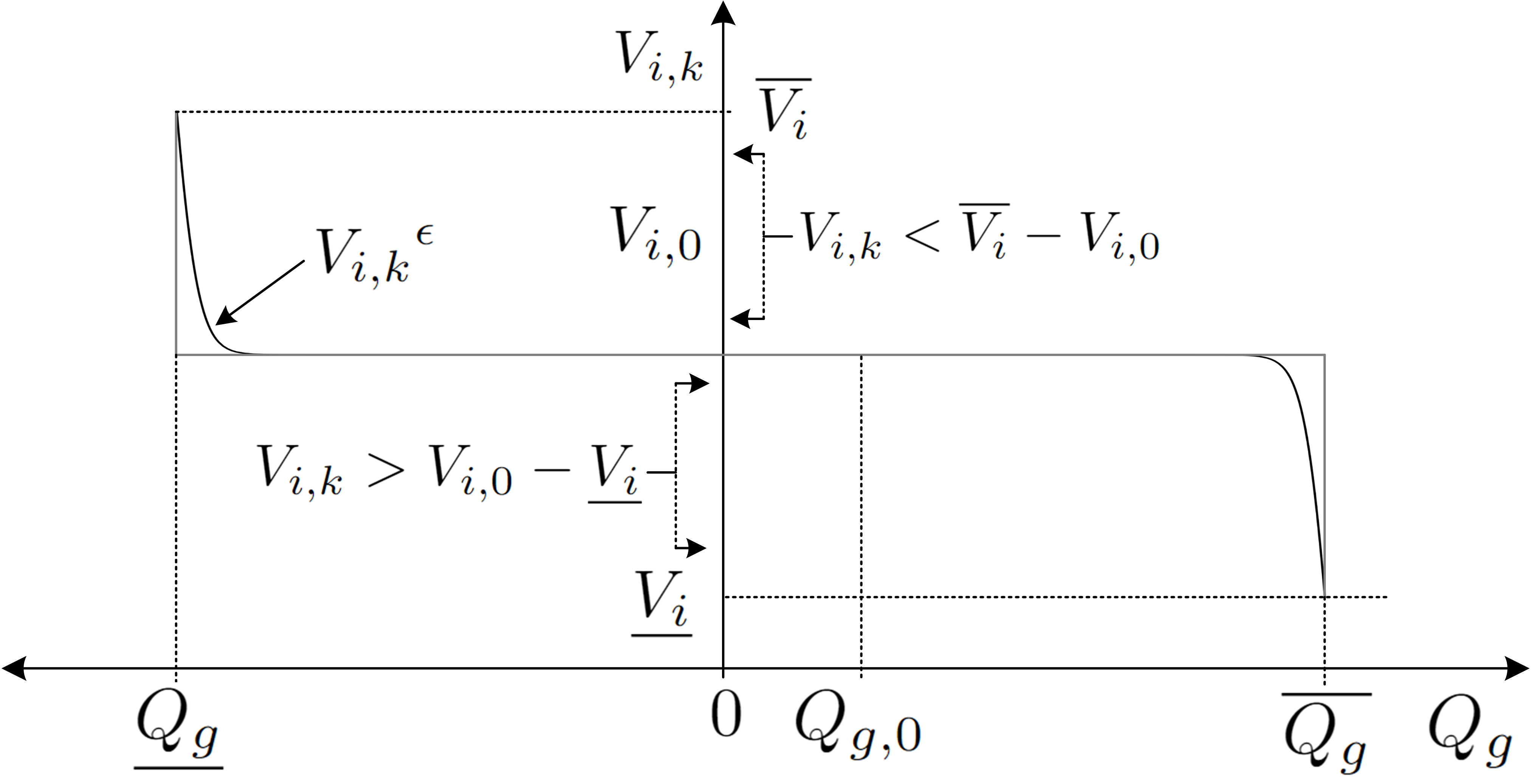}
    \caption{Original and approximate generator reactive power response functions.}
    \label{fig:genQ}
\end{figure}

\section{Proposed Smoothened Droop Functions}
 The PF, OPF, and SCOPF models can be accessed from our previous works \cite{PF}, \cite{pmacdc}, \cite{SCOPF} respectively. Only, the generator PQ and converter droop control constraints with their proposed smooth formulations are presented in this section. 

\subsection{Smooth Approximation of Piece-Wise Linear Constraints}
Consider the ReLU function $f(x) = \max\{0,x\}$, such that $x \in \mathbb{R}$, that is not differentiable at $x = 0$. For any positive value $\epsilon > 0$, $f(x)$ can be approximated by the smooth softplus function \cite{056} with sharpness setting $\epsilon$,
\begin{equation}
    \label{eq:smooth_fcn}
    f^\epsilon(x) := \epsilon \ln(1 + \exp^{x/\epsilon}) \text{ for } x \in \mathbb{R},
\end{equation}
with verifiable bounds on approximation error,
\begin{equation}
    \label{eq:smooth_fcn_error}
    f^\epsilon(x)-\epsilon \ln 2 \leq f(x) \leq f^\epsilon(x).
\end{equation}


\subsection{Generator PQ Control in SCOPF Model}
A generator responding in a contingency event $k$ adjusts its active power output by a common perturbation denoted by variable $\Delta_k$ according to its predefined (offline) participation factor $\alpha_g$ (droop slope) until it hits an operational bound $\Pgmin$, or $\Pgmax$ as shown in Fig. \ref{fig:genP}. This is referred to as the generator's post-contingency frequency response. The active power output $P_{g,k}$ of a responding generator $g$ in contingency $k$ is given as,
\begin{multline} 
    \label{eq38} 
    P_{g,k}\!=\!
    \begin{cases} 
        \Pgmax, & P_{g,k}\!\geq\!\Pgmax, \\ 
        P_{g,0}\!+\!\alpha_g \Delta_k, & \Pgmin\!\leq\!P_{g,k}\!\leq\!\Pgmax, \\
        \Pgmin, & P_{g,k}\!\leq\!\Pgmin, \\ 
    \end{cases} 
    \forall g \in \mathcal{G}, k \in \setK.
\end{multline}
A generator that is online in a contingency but is not selected to respond to that contingency maintains its active power output from the base case,
\begin{equation} 
    \label{eq39}
    {P^{\text{g}}_{n,k} = P_{g,0}, \;\; \forall g \in \mathcal{G}_k \notin \mathcal{G}_k^{\text{nr}}, \quad k \in \setK}.
\end{equation}
In a contingency event $k$, a responding generator $g$ tries to maintain its voltage magnitude $V_{i,k}$ equal to the base case voltage $V_{i,0}$ by adjusting its reactive power $Q_{g,k}$ unless it hits operational bounds $\underline{Q_{g}}$ or $\overline{Q_{g}}$ as shown in Fig. \ref{fig:genQ}. This is referred to as PV/PQ bus switching control and modeled,
\begin{flalign} 
    \label{eq40}
    \begin{cases}
    V_{i,k} \leq V_{i,0}, & Q_{g,k} = \overline{Q_{g}}, \\
    V_{i,k} = V_{i,0}, & \underline{Q_{g}} \leq Q_{g,k} \leq \overline{Q_{g}},\; \forall gi \in \mathcal{G} \times \mathcal{I}, k \in \setK.\\
    V_{i,k} \geq V_{i,0}, & Q_{g,k} = \underline{Q_{g}}, \\
    \end{cases}
\end{flalign}
We convert conditional constraints (\ref{eq38}) and (\ref{eq40}) into continuously differentiable functions using smooth approximation \eqref{eq:smooth_fcn}. The constraint (\ref{eq38}) can be written,
\begin{multline} 
    \label{eq41}
    P_{g,k} = \max\{\underline{P_{g}}, -\max\{-\overline{P_{g}}, -P_{g,0} -\alpha_g \Delta_k \} \},\\
    \forall g \in \mathcal{G}, k \in \setK. 
\end{multline}
By applying \eqref{eq:smooth_fcn}, constraint (\ref{eq41}) can be expressed,
\begin{multline} 
    \label{eq42}
    P_{g,k}^{ \epsilon} := \underline{P_{g}} +
    \epsilon \ln \left( 1 + \frac{\exp^{\big((\overline{P_{g}} - \underline{P_{g}})/\epsilon\big)}}{1 + \exp^{\big((\overline{P_{g}} - P_{g,0} - \alpha_g \Delta_k)/\epsilon\big)}} \right),\\
    \forall g \in \mathcal{G}, k \in \setK. 
\end{multline}
The parameter $\epsilon >0$ can be used to regulate the degree of approximation accuracy as shown in Fig. \ref{fig:genP}, the approximation gap narrows for $\epsilon_1 < \epsilon_2$ and further as $\epsilon \rightarrow 0$.
Similarly, with generator-node mapping $\mathcal{T}^g \subset \mathcal{G} \times \mathcal{I}$,  (\ref{eq40}) can be stated,
\begin{multline} \label{eq43}
    V_{i,k}^{\epsilon} := V_{i,0} +  \max \{ V_{i,k}^+ -{Q_{g,k}} + \underline{Q_g}, 0 \} \\ -  \max \{ {V_{i,k}^- +Q_{g,k}} - {\overline{Q_g}}, 0 \}, \; \forall gi \in \mathcal{T}^g , k \in \setK,
\end{multline}
where $V_{i,k}^+$ and $V_{i,k}^-$ represent the upper and lower margins of base case voltage $V_{i,0}$ from the bounds as,
\begin{align}
    0 &\leq V_{i,k}^+ \leq \overline{V_i} - V_{i,0}, \;\; \forall gi \in \mathcal{T}^g, k \in \setK, \label{eq44}\\
    0 &\leq V_{i,k}^- \leq V_{i,0} - \underline{V_i}, \;\; \forall gi \in\mathcal{T}^g, k \in \setK. \label{eq45}
\end{align}
The necessary conditions for (\ref{eq43}) are the voltage magnitude $(\underline{V_i}, \overline{V_i})$ and generator reactive power bounds at these buses.
By applying \eqref{eq:smooth_fcn}, the constraint \eqref{eq43} can be expressed as,
\begin{multline}\label{eq48}
V_{i,k}^{\epsilon} = V_{i,0} + \epsilon \ln \bigg(1 + \exp^{\big((V_{i,k}^+ - Q_{g,k} + \underline{Q_{g}} )/\epsilon\big)} \bigg) \\
- \epsilon \ln \bigg(1 + \exp^{\big((V_{i,k}^- + Q_{g,k} - \overline{Q_{g}} )/\epsilon\big)}\bigg),
\scalemath{0.79}{\forall gi \in \mathcal{G} \times \mathcal{I}, k \in \setK.}
\end{multline}

\begin{figure}[!t]
    \centering
    \includegraphics[width=2.8in]{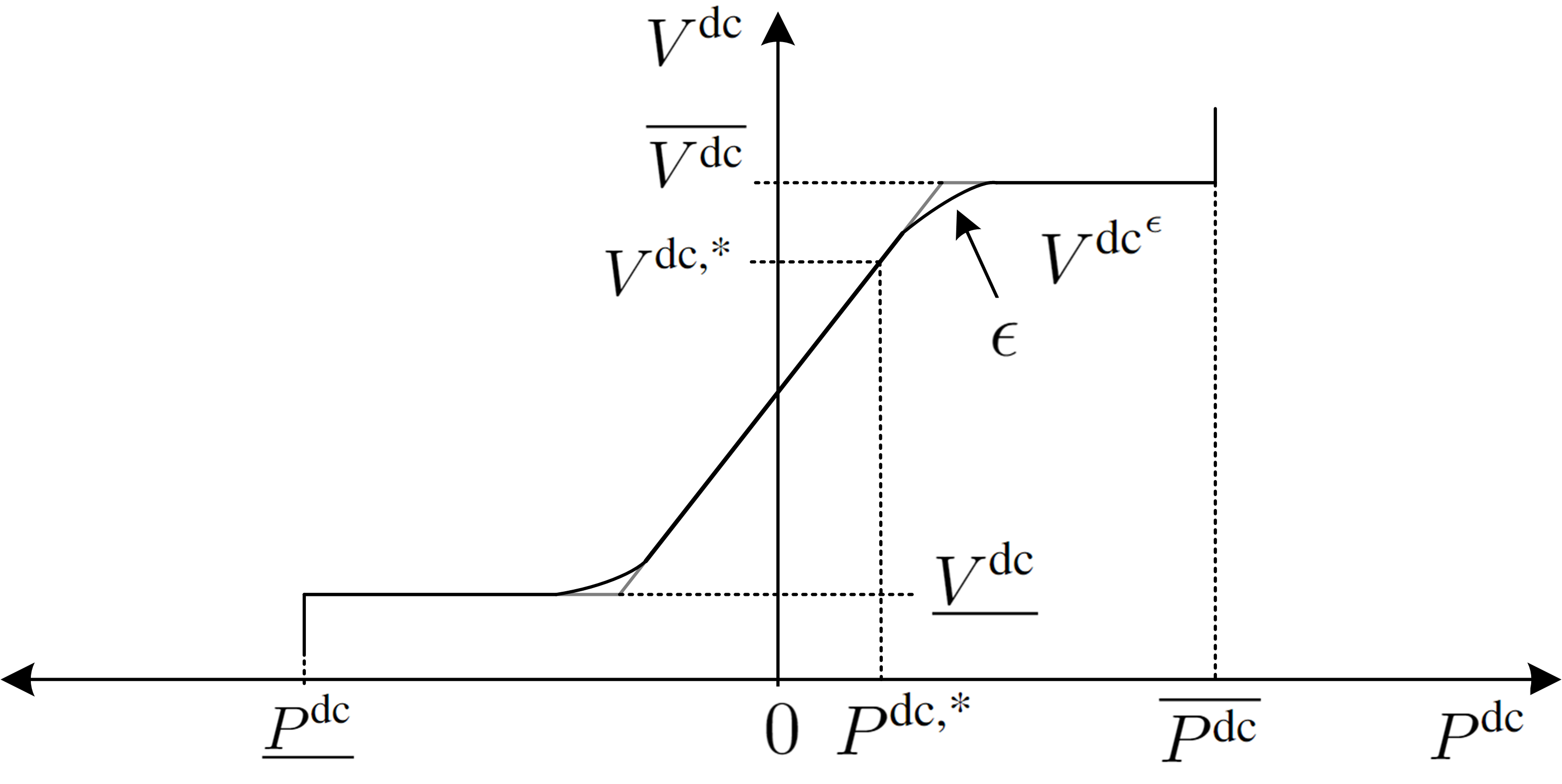}
    \caption{Original and approximate P--Vdc linear droop functions.}
    \label{fig:PVdcdroop}
\end{figure}
\begin{figure}[!t]
    \centering
    \includegraphics[width=2.75in]{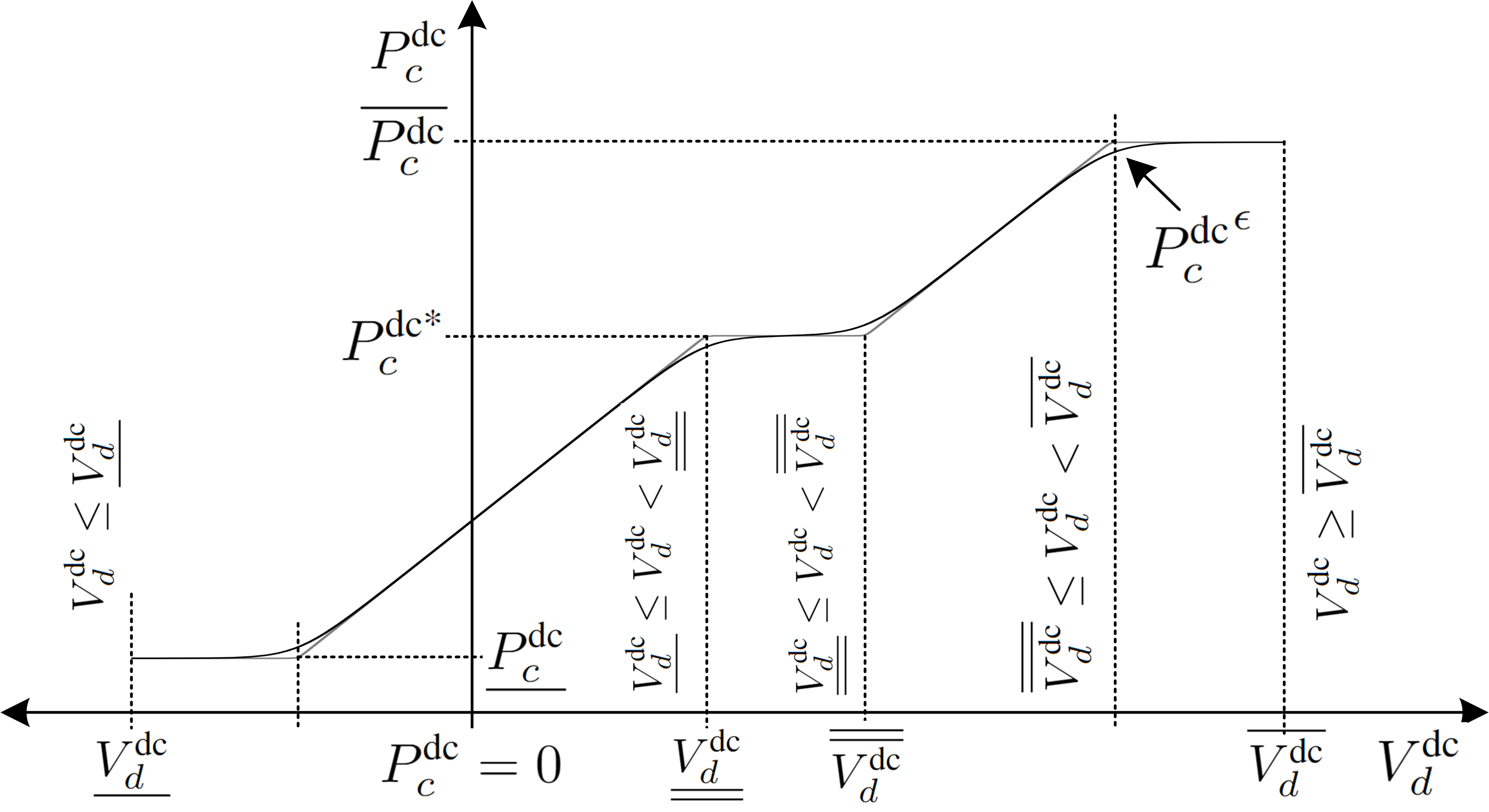}
    \caption{Original and approximate P--Vdc piecewise linear droop functions.}
    \label{fig:PVdcPWdroop}
\end{figure}

\subsection{HVDC Converter P--Vdc Droop Control Model}
The P--Vdc linear droop control as shown in Fig. \ref{fig:PVdcdroop} is expressed by a relation between dc voltage $V^{\text{dc}}_d$ at dc grid bus $d \in \mathcal{D}$ and the dc side power $P^{\text{dc}}_{c}$ of converter $c \in  \mathcal{C}$, with the mapping $\mathcal{T}^{c} \subset \mathcal{C} \times \mathcal{D}$ and $ \forall cd \in \mathcal{T}^{c} $,
\begin{multline} 
    \label{eq:linearDroop} 
    P^{\text{dc}}_{c}\!=\!
    \begin{cases} 
        P^{\text{dc*}}_{c}\!-\!\frac{\zeta(P^{\text{dc*}}_{c})} {k^{\text{droop}}_c} (\overline{V^{\text{dc}}_d}-V_{d}^{\text{dc*}}), & V^{\text{dc}}_d\!\geq\!\overline{V^{\text{dc}}_d}, \\ 
        P^{\text{dc*}}_{c}\!-\!\frac{\zeta(P^{\text{dc*}}_{c})} {k^{\text{droop}}_c} (V_{d}^{\text{dc}}-V_{d}^{\text{dc*}}), & \underline{V^{\text{dc}}_d}\!\leq\!V_{d}\!\leq\!\overline{V^{\text{dc}}_d}, \\
        P^{\text{dc*}}_{c}\!-\!\frac{\zeta(P^{\text{dc*}}_{c})} {k^{\text{droop}}_c} (\underline{V^{\text{dc}}_d}-V_{d}^{\text{dc*}}), & V^{\text{dc}}_d\!\leq\!\underline{V^{\text{dc}}_d}, \\ 
    \end{cases}
\end{multline}
where $k^{\text{droop}}_c$ is the droop coefficient, $ P^{\text{dc*}}_c$ denotes the reference dc side power, and $V^{\text{dc*}}_d$ is the reference dc voltage. The $\zeta$ function returns the sign of the reference dc side power to represent both inverter and rectifier modes of operation.

The constraint \eqref{eq:linearDroop} can be rewritten in the form of max functions explicitly for either of the sign conventions based on a given reference value, 
\begin{multline} 
    \label{eq:linearDroopMax}
P^{\text{dc}}_{c} =
    \begin{cases} 
    \max \big\{ \underline x, \max \big\{\overline x, P^{\text{dc}}_{c}(V^{\text{dc}}_d) \big\}\big\},  &P^{\text{dc*}}_{c} > 0, \\
    \max \big\{ \underline x, -\max \big\{-\overline x, - P^{\text{dc}}_{c}(V^{\text{dc}}_d) \big\}\big\},  &P^{\text{dc*}}_{c} < 0, \\
    \end{cases}\\
    \forall cd \in \mathcal{T}^{c}.
\end{multline}

Here, the converter dc power bounds are expressed as $ \underline x = \max \{P^{\text{dc}}_{c}(\underline{V^{\text{dc}}_d}), \underline{P^{\text{dc}}_{c}} \} $, and $\overline x = \min\{P^{\text{dc}}_{c}(\overline{V^{\text{dc}}_d}), \overline{P^{\text{dc}}_{c}} \}$ using \eqref{eq:linearDroop}.

By applying \eqref{eq:smooth_fcn}, constraint (\ref{eq:linearDroopMax}) can be expressed,
\begin{multline}
    \label{eq:smoothlinearDroop}
    {P^{\text{dc}}_{c}}^\epsilon :=
    \begin{cases}
    \underline{x} - \epsilon \ln \bigg( 1 + \frac{\exp^{((-\overline{x} + \underline{x})/\epsilon)}}{1 + \exp^{((-\overline{x} + P^{\text{dc}}_{c}(V^{\text{dc}}_d))/\epsilon)}} \bigg), &P^{\text{dc*}}_{c} > 0,\\
    \underline{x} + \epsilon \ln \bigg( 1 + \frac{\exp^{((\overline{x} - \underline{x})/\epsilon)}}{1 + \exp^{((\overline{x} - P^{\text{dc}}_{c}(V^{\text{dc}}_d))/\epsilon)}} \bigg), &P^{\text{dc*}}_{c} < 0,\\
    \end{cases}\\
    \forall cd \in \mathcal{T}^{c}.
\end{multline}

The P--Vdc piecewise linear droop control with power deadband and voltage limits as shown in Fig. \ref{fig:PVdcPWdroop} can be defined,
\begin{multline} 
    \label{eq:PWlinearDroop}
    P^{\text{dc}}_{c} = P^{\text{dc*}}_{c} +
    \begin{cases}
        \frac{(\overline{V^{\text{dc}}_d} - V^{\text{dc}}_d)}{{k^{\text{droop}}_c}}    + \frac{(\overline{\overline{V^{\text{dc}}_d}} - \overline{V^{\text{dc}}_d})}{k^{\text{droop}}_c} ,
            & V^{\text{dc}}_d \geq \overline{V^{\text{dc}}_d}, 
        \\
        \frac{(\overline{\overline{V^{\text{dc}}_d}} - V^{\text{dc}}_d)}{k^{\text{droop}}_c} , 
            &  \overline{\overline{V^{\text{dc}}_d}} \leq V^{\text{dc}}_d < \overline{V^{\text{dc}}_d} , 
        \\
        0, 
            &  \underline{\underline{V^{\text{dc}}_d}} \leq V^{\text{dc}}_d < \overline{\overline{V^{\text{dc}}_d}},
        \\
        \frac{(\underline{\underline{V^{\text{dc}}_d}} - V^{\text{dc}}_d)}{k^{\text{droop}}_c} , 
            &  \underline{V^{\text{dc}}_d} \leq V^{\text{dc}}_d < \underline{\underline{V^{\text{dc}}_d}}, 
        \\
        \frac{(\underline{V^{\text{dc}}_d} - V^{\text{dc}}_d)}{{k^{\text{droop}}_c}}  +  \frac{(\underline{\underline{V^{\text{dc}}_d}} - \underline{V^{\text{dc}}_d} )}{k^{\text{droop}}_c} , 
            &  V^{\text{dc}}_d \leq \underline{V^{\text{dc}}_d}, \\
    \end{cases}
    \\
    \forall cd \in \mathcal{T}^{c}.
\end{multline}

The converter operates in a constant power mode when the dc voltage is between $\underline{\underline{V^{\text{dc}}_d}}$ and $\overline{\overline{V^{\text{dc}}_d}}$. When the voltage exceeds these limits, the droop control adjusts the power output of the converter until the voltage hits its bounds ${\underline{V^{\text{dc}}_d}}$ and ${\overline{V^{\text{dc}}_d}}$. Similar to \eqref{eq:linearDroop}, a smooth function for \eqref{eq:PWlinearDroop} with voltage limits is,
\begin{multline}
\label{eq:PWsmooth}
    {P^{\text{dc}}_{c}}^\epsilon := -\zeta(P^{\text{dc*}}_{c}) P^{\text{dc*}}_{c} + 1/k^{\text{droop}}_c \bigg({\overline{V^{\text{dc}}_d} + \underline{V^{\text{dc}}_d} - \underline{\underline{V^{\text{dc}}_d}} - \overline{\overline{V^{\text{dc}}_d}}} \bigg)  
    \\  
    + \epsilon \ln \bigg(1 + \exp^{\big((1/k^{\text{droop}}_c {(\overline{\overline{V^{\text{dc}}_d}} - V^{\text{dc}}_d)} - \overline{V^{\text{dc}}_d} + V^{\text{dc}}_d)/\epsilon\big)} \bigg) 
    \\
    -\epsilon \ln \bigg(1 + \exp^{\big((1/k^{\text{droop}}_c (\overline{V^{\text{dc}}_d} - V^{\text{dc}}_d) - 2\overline{V^{\text{dc}}_d} + \overline{\overline{V^{\text{dc}}_d}} + V^{\text{dc}}_d)/\epsilon\big)}\bigg) 
    \\
    -\epsilon \ln \bigg(1 + \exp^{\big( (- 1/k^{\text{droop}}_c (\underline{\underline{V^{\text{dc}}_d}} - V^{\text{dc}}_d) + \underline{V^{\text{dc}}_d} - V^{\text{dc}}_d) / \epsilon \big) } \bigg) 
    \\
    +\epsilon \ln \bigg(1 + \exp^{\big( (- 1/k^{\text{droop}}_c (\underline{V^{\text{dc}}_d} - V^{\text{dc}}_d) + 2 \underline{V^{\text{dc}}_d} - \underline{\underline{V^{\text{dc}}_d}} - V^{\text{dc}}_d) /\epsilon\big)} \bigg), 
    \\
    P^{\text{dc*}}_{c} < 0, \forall cd \in \mathcal{T}^{c}.
\end{multline}
When $P^{\text{dc*}}_{c} > 0$, $-{P^{\text{dc}}_{c}}^\epsilon (V^{\text{dc}}_d)$ can be implimented. For the PF application, \eqref{eq:linearDroop} and \eqref{eq:PWlinearDroop} are implemented without limits and thus, this five-step piecewise linear droop function is extended beyond voltage bounds. Furthermore, power limits can also be imposed on \eqref{eq:PWsmooth} similar to \eqref{eq:linearDroopMax} and \eqref{eq:smoothlinearDroop}. Finally, the same smooth functions can be adapted for Q--Vac droop control.


\section{Simulation Results}
Two test datasets, \texttt{case5} and \texttt{case67} are used to perform PF, OPF, and, SCOPF simulations. The \textsc{Ipopt} and \textsc{Juniper} solvers are used for optimization models on a PC with Intel Core i9-11950H with 64 GB  RAM. The test datasets, PF, OPF, and SCOPF models are publicly available as an open-source library PowerModelsACDCsecurityconstrained.jl’\footnote{\scriptsize{https://github.com/csiro-energy-systems/PowerModelsACDCsecurityconstrained.jl}}. 

\subsection{Comparing Proposed Models with/out Integer Variables}
The section compares the proposed smooth approximation-based NLP models with the MINLP models for the PF, OPF and SCOPF problems considering steady-state voltage regulation in ac-dc power systems using P--Vdc and Q--Vac droop controls and generator PQ response functions. Five scenarios with the following response functions are studied:
\begin{enumerate}[i)]
    \item generator PQ droop;
    \item generator PQ, and converter P-Vdc droops;
    \item generator PQ, and converter P-Vdc w. deadband droops;
    \item generator PQ, and converter P-Vdc, Q-Vac droops;
    \item generator PQ, and converter P-Vdc, Q-Vac w. deadband droops;
\end{enumerate}

The numerical results, as presented in Table \ref{table:comparisonPfOpfScopf}, illustrate that the proposed smooth approximation-based NLP models significantly outperform the MINLP models in computational efficiency. While the objective function values remain largely consistent across the models, the incorporation of piecewise linear droop curves with deadbands in the MINLP models introduces additional binary variables (bvar), which increases the computation time. Notably, scenario v, characterized by the highest number of binary variables, exhibits the greatest computational burden. In the context of the SCOPF problem, the addition of each contingency substantially increases the number of continuous (cvar) and binary variables, rendering the problem increasingly challenging to solve. In the 67-bus case, in scenarios iii, iv, and v, MINLP SCOPF models could not be solved due to limited computational resources, whereas, the proposed NLP models efficiently solve these scenarios.

\begin{table*} 
\centering
\caption{Comparison of Proposed Smooth Approximation-based PF, OPF, and SCOPF NLP Models with MINLP OPF and SCOPF Models}
\label{table:comparisonPfOpfScopf}
\begin{adjustbox}{width=0.99\linewidth,center}
\begin{tabular}{l c c c c c c c c c c |c| c c c c c c c c c}
\hline 
case & sc. && \multicolumn{4}{c}{OPF MINLP Model} && \multicolumn{2}{c}{OPF Propsoed NLP Model}&&  PF  && \multicolumn{4}{c}{SCOPF MINLP Model}  && \multicolumn{3}{c}{SCOPF Propsoed NLP Model} \\
 \cline{4-7} \cline{9-10} \cline{12-12} \cline{14-17} \cline{19-21}
&                        && obj ($\$$)   & time (s) & $\#$ bvar & $\#$ cvar && obj ($\$$)   & time (s)         && time (s) && obj ($\$$)   & time (s)  & $\#$ bvar & $\#$ cvar && obj ($\$$)   & time (s)  & $\#$ cvar 
\\ \hline
\multirow{5}{*}{5} & i   && 30        & 10.0  & 25  & 147 && 30        & 0.13  && --  && 53.5276   &244.4    & 88  & 1880  && 52.4455   & 13.644 & 1694   \\
                   & ii  && 30        & 16.4  & 45  & 147 && 30        & 0.10  && 0.11 && 53.6021   &296.7    & 126 & 1880  && 53.6628   & 10.445 & 1694   \\
                   & iii && 30        & 27.2  & 55  & 147 && 30        & 2.14  && 0.13 && 54.0348   &288.4    & 153 & 1880  && 53.159    & 10.824 & 1694   \\
                   & iv  && 30        & 21.4  & 65  & 147 && 30        & 0.16  && 0.13 && 33853.1   &520.2    & 198 & 1880  && 32470.9   & 9.961  & 1694   \\
                   & v   && 30        & 35.2  & 85  & 147 && 30        & 3.40  && 0.15 && 96558.4   &578.5    & 218 & 1880  && 96487.7   & 18.732 & 1694   \\ \hline 
\multirow{5}{*}{67}& i   && 1.233e5 & 162.5  & 45  & 742 && 1.234e5 & 0.57  && --  && 1.246e5 &13178.4    & 512 & 6366  && 1.257e5 & 123.44 & 26756  \\
                   & ii  && 1.234e5 & 215.9  & 72  & 742 && 1.234e5 & 0.80  && 0.23 && 1.238e5 &14376.3    & 768 & 9549  && 1.345e5 & 145.62 & 26756  \\
                   & iii && 1.245e5 & 286.4  & 90  & 742 && 1.234e5 & 0.73  && 0.21 &&  ---      & ---     & --- & ---   && 1.246e6 & 106.33 & 21353  \\
                   & iv  && 1.247e5 & 1555.3 & 99  & 742 && 1.245e5 & 1.48  && 0.20 && ---       & ---     & --- & ---   && 1.340e5 & 126.81 & 29429  \\
                   & v   && 1.246e5 & 1836.1 & 135 & 742 && 1.245e5 & 1.14  && 0.22 && ---       & ---     &---  & ---   && 1.340e7 & 116.50 & 26756  \\ \hline 
\end{tabular}
\end{adjustbox}
\end{table*}

\begin{figure}[!t]
    \centering
    \includegraphics[height=1.9in, width=2.7in]{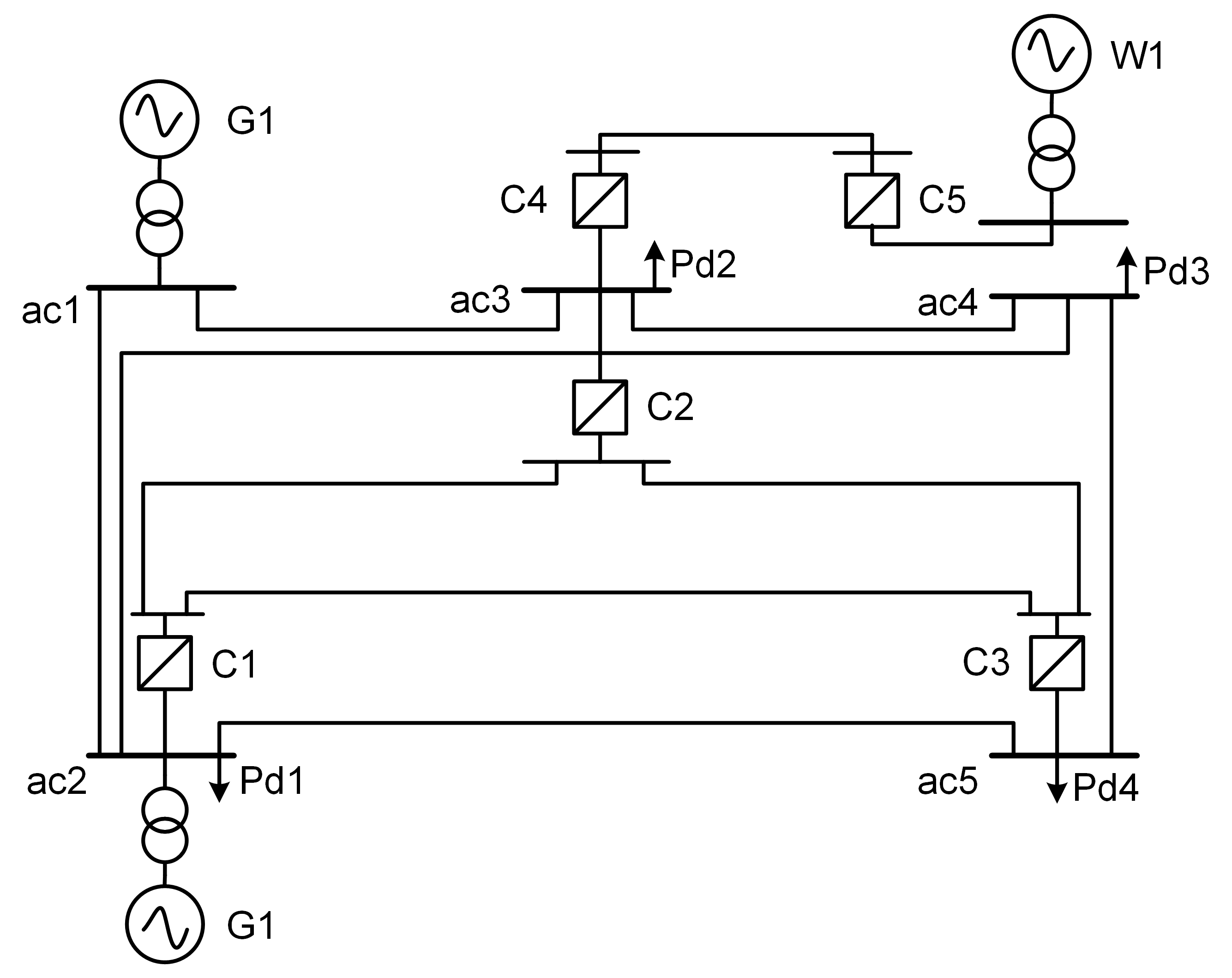}
    \caption{\texttt{case5} ac-dc test network with HVDC connected offshore wind farm.}
    \label{fig:case5}
\end{figure}

\begin{figure*}[ht!]
    \centering
    \includegraphics[width=1\linewidth]{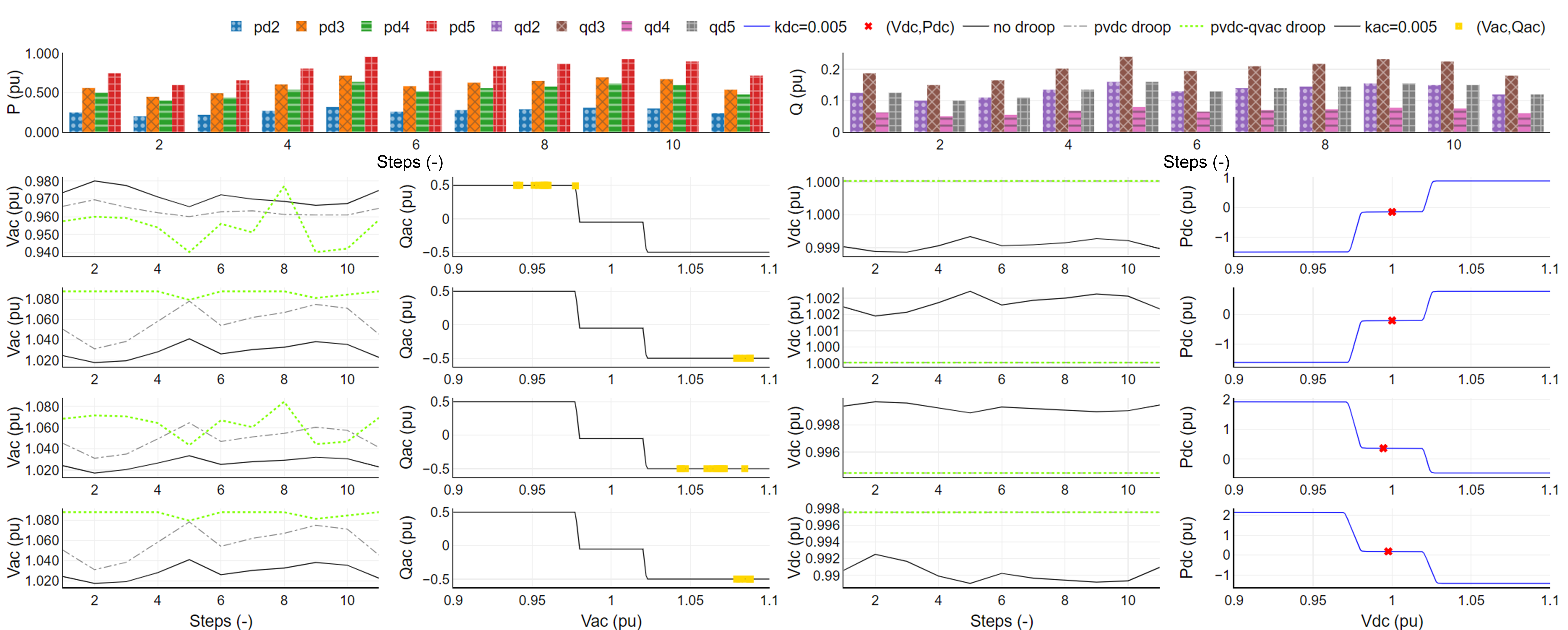}
    \caption{Voltage regulation with and without P--Vdc only and P--Vdc, Q--Vac droop control for \texttt{case5} OPF at different load variations.}
    \label{fig:opfdroop}
\end{figure*}

\subsection{Voltage Regulation with P--Vdc and Q--Vac droop control}

The proposed smooth approximation-based OPF models with and without P--Vdc only and P--Vdc, Q--Vac droop control are solved on \texttt{case5} ac-dc power system with different load variations and the results are presented in Fig. \ref{fig:opfdroop}. These results indicate that the implementation of P--Vdc and Q--Vac droop control substantially enhances voltage regulation. Furthermore, the findings underscore that the limited reactive power capacity of the HVDC converter necessitates allowing voltage drops only when the voltage reaches its bounds. This constraint highlights the importance of effectively managing reactive power within the system to maintain voltage stability under varying load conditions.   

\begin{figure*}[t]
    \centering
    \includegraphics[width=1\linewidth, height=6.4cm]{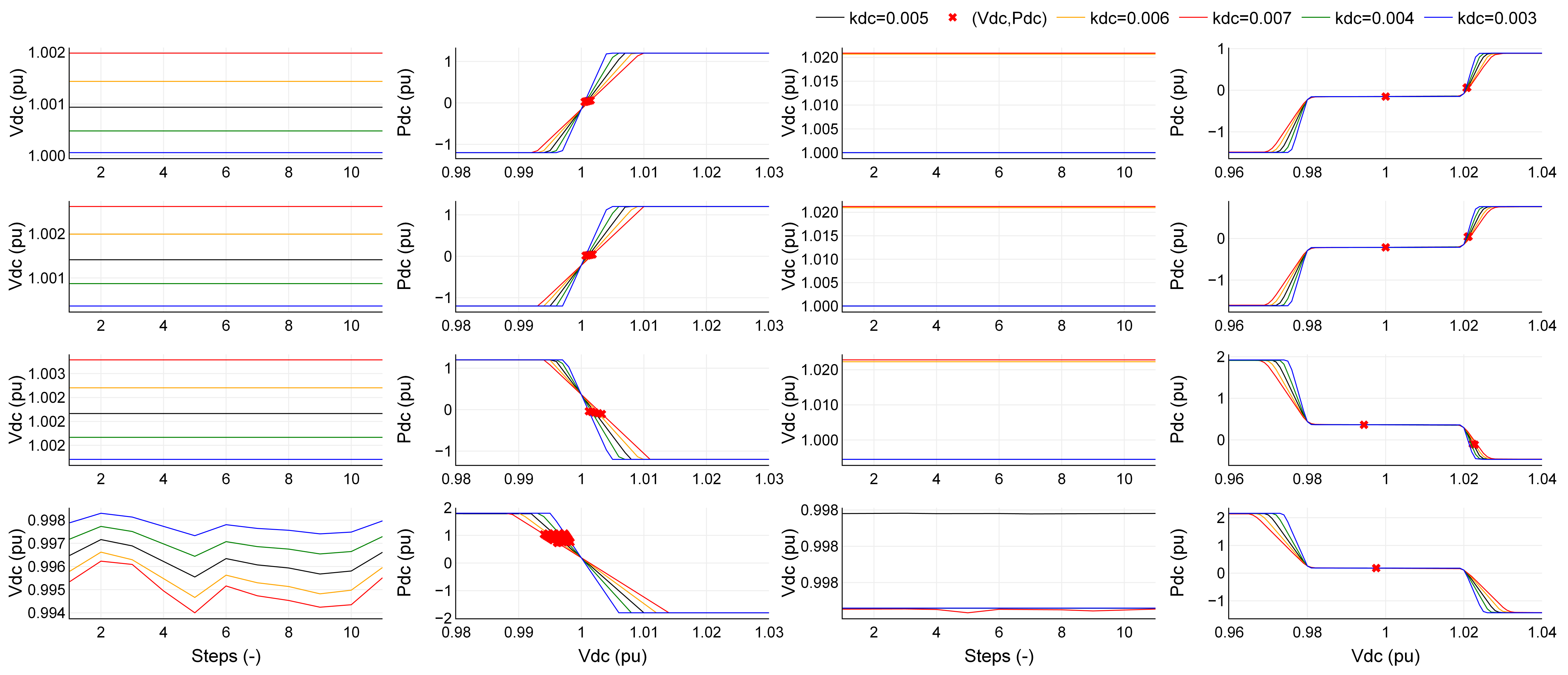}
    \caption{HVDC converters P--Vdc droop control with and without deadband and different slope coefficients for \texttt{case5} OPF at different load variations.}
    \label{fig:kdroop}
\end{figure*}

\subsection{Impact of Droop Coefficient and Deadband}

This section evaluates the impact of different P--Vdc droop slopes and the implementation of droop control with and without deadband on voltage regulation under varying load conditions in the \texttt{case5} ac-dc power system as shown in Fig. \ref{fig:kdroop}. The droop control without a deadband employs a fixed droop coefficient, balancing the trade-off between voltage regulation stiffness and power-sharing accuracy. A small droop coefficient such as 0.003 in comparison to 0.007 improves voltage regulation but sacrifice power-sharing and vice versa. Typically, system operators allow a maximum 10\% DC voltage variation, allocating 5\% for the line drops, thus limiting the voltage range available for droop control and consequently the droop coefficient. However, the droop control without a deadband cannot operate the converter in constant power control. In contrast, droop control with a deadband maintains the converters in constant power control when the voltage is between 0.98 and 1.02 pu, therefore providing better power-sharing under higher loads and better voltage regulation under light loads.

\begin{figure*}[t]
    \centering
    \includegraphics[width=1\linewidth]{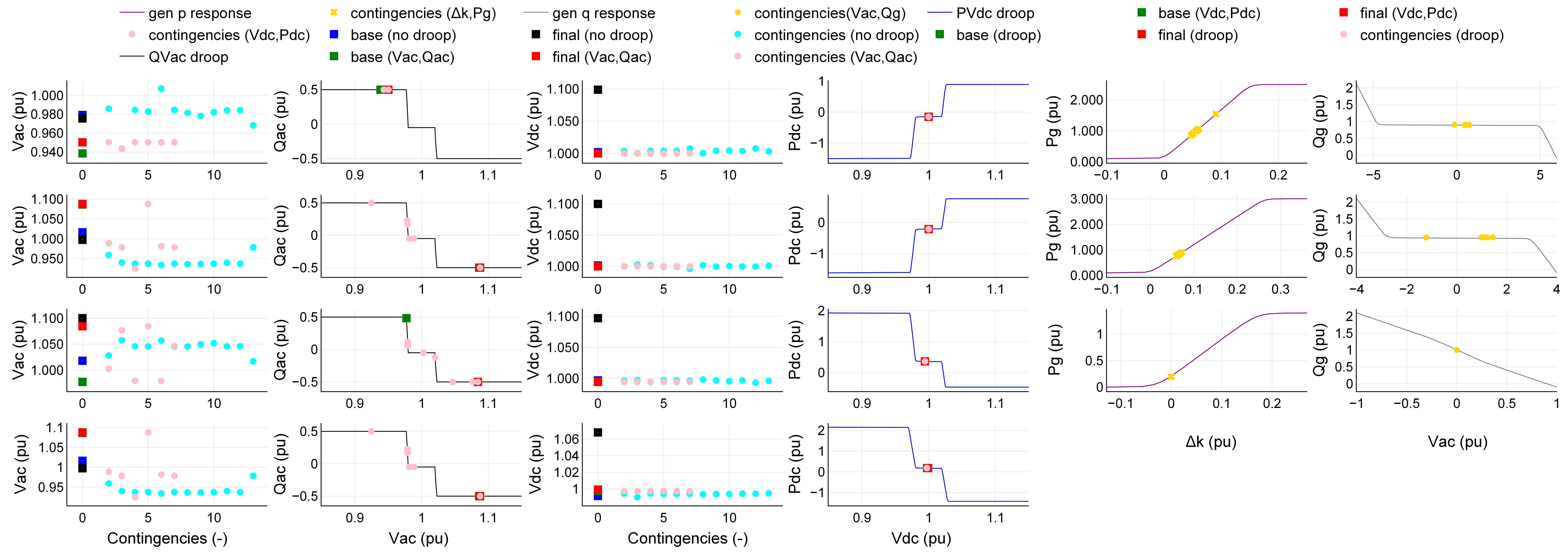}
    \caption{Voltage regulation with and without P--Vdc, Q--Vac droop control and generator PQ response for \texttt{case5} SCOPF contingency events.}
    \label{fig:droopScopf}
\end{figure*}

\subsection{Voltage Regulation with P--Vdc, Q--Vac Droop control and Generator PQ Response in Contingency Events}
The proposed smooth approximation allows efficient implementation of the P--Vdc, Q--Vac droop control, and generator PQ response functions in the SCOPF problem. As illustrated in Fig. \ref{fig:droopScopf}, the P--Vdc, Q--Vac droop control for \texttt{case5} ac-dc power system efficiently regulates voltage and manages power-sharing during various contingency events such as generator, ac transmission line, transformer, converter, and dc transmission line contingencies. Specifically generator 3, as shown in Fig. \ref{fig:droopScopf}, is an offshore wind farm connected to the mainland via a HVDC system with no reactive power capability. The P--Vdc droop control ensures the converters operate in constant power control mode by leveraging the sufficient active power response from generators during contingencies, thus maintaining the base-case voltage throughout the events and achieving the final solution. Meanwhile, the Q--Vac droop control considers converters' limited reactive power capacity, effectively regulating ac bus voltages and adhering to reactive power and voltage constraints.
\vspace{-0.3cm}
\section{Conclusions}
The proposed smooth approximation-based modeling enables implementation of generator PQ response, converters P--Vdc, Q--Vac droop control including voltage limits, power deadbands, and power limits in PF, OPF, and SCOPF applications. These models employ smooth (continuous and differentiable) functions that can be efficiently solved using Newton’s method or derivative-based interior point solvers. Numerical results demonstrate a substantial reduction in computational time compared to existing MINLP models. Furthermore, these models are fully parameterized and tuneable, offering flexibility for enhanced system efficiency, while the findings highlight improved voltage regulation and power-sharing during load variations and contingency events.

\vspace{0cm}
\footnotesize{
\bibliographystyle{IEEEtranN}
\bibliography{References}
}


\end{document}